\documentstyle[a4,11pt]{article}
\begin{document}
\begin{titlepage}
\vskip 2cm
\begin{flushright}
Preprint CNLP-1994-02
\end{flushright}
\vskip 2cm
\begin{center}
{\bf
SOLITON EQUATIONS IN 2+1 DIMENSIONS
AND DIFFERENTIAL GEOMETRY OF CURVES/SURFACES}
\footnote{Preprint
CNLP-1994-02. Alma-Ata. 1994 }
\end{center}
\vskip 2cm
\begin{center}
Ratbay MYRZAKULOV
\footnote{E-mail: cnlpmyra@satsun.sci.kz}
\end{center}
\vskip 1cm

\begin{center}
 Centre for Nonlinear Problems, PO Box 30, 480035, Alma-Ata-35, Kazakhstan
\end{center}

\begin{abstract}
Some aspects of the relation between differential geometry of curves and
surfaces  and (2+1)-dimensional soliton equations are discussed.
For the (2+1)-dimensional case,
self-cordination of geometrical formalism with the Hirota's bilinear
method is established.
 A connection
between supersymmetry, geometry and soliton equations is also considered.
\end{abstract}


\end{titlepage}

\setcounter{page}{1}
\newpage

\tableofcontents
\section{Introduction}
Consider the curve in 3-dimensional space. Equations of such curves,
following [1] we can write in the form
$$
\left ( \begin{array}{ccc}
{\bf e}_{1} \\
{\bf e}_{2} \\
{\bf e}_{3}
\end{array} \right)_{x}= C
\left ( \begin{array}{ccc}
{\bf e}_{1} \\
{\bf e}_{2} \\
{\bf e}_{3}
\end{array} \right)
\eqno(1a)
$$
$$
\left ( \begin{array}{ccc}
{\bf e}_{1} \\
{\bf e}_{2} \\
{\bf e}_{3}
\end{array} \right)_{t}= G
\left ( \begin{array}{ccc}
{\bf e}_{1} \\
{\bf e}_{2} \\
{\bf e}_{3}
\end{array} \right) \eqno(1b)
$$
where
$$
C =
\left ( \begin{array}{ccc}
0             & k     &  0 \\
-\beta k      & 0     & \tau  \\
0             & -\tau & 0
\end{array} \right) ,\quad
G =
\left ( \begin{array}{ccc}
0       & \omega_{3}  & -\omega_{2} \\
-\beta\omega_{3} & 0      & \omega_{1} \\
\beta\omega_{2}  & -\omega_{1} & 0
\end{array} \right) \eqno(2)
$$
Here
$$
{\bf e}_{1}^{2}=\beta = \pm 1, {\bf e}_{2}^{2}={\bf e}^{2}_{3}=1 \eqno(3)
$$
Note that equation (1a) is the Serret - Frenet equation (SFE). So we have
$$
C_t - G_x + [C, G] = 0 \eqno(4)
$$
or
$$
k_{t} - \omega_{3x}  - \tau \omega_{2} = 0 \eqno (5a)
$$
$$
\omega_{2x} - \tau \omega_{3} + k \omega_{1} = 0 \eqno (5b)
$$
$$
\tau_{t} - \omega_{1x} + \beta k \omega_{2}  = 0. \eqno (5c)
$$

We now consider the isotropic Landau-Lifshitz equation (LLE)
$$
{\bf S}_{t} = {\bf S}\wedge {\bf S}_{xx}. \eqno(6)
$$
If
$$
{\bf e}_{1} \equiv {\bf S}  \eqno(7)
$$
then
$$
q=\frac{k}{2}e^{i\partial^{-1}_{x} \tau} \eqno(8)
$$
satisfies the NLSE
$$
iq_{t}+q_{xx}+2\beta \mid q\mid^{2}q =0. \eqno(9)
$$
This equivalence between the LLE (6) and the NLSE (9) we call the Lakshmanan
equivalence or L-equivalence [1]. These results for the case
$\beta =+1$ was obtained in [2] and for the case $\beta =-1$ in [1].
Note that between these equations also take places gauge equivalence (G-equivalence) [6].

In this paper, starting from Lakshmana's idea,
we will discuss some aspects of the relation between
differential geometry of curves and surfaces and (2+1)-dimensional soliton
equations. Before this, in [1] we proposed some approaches to this problem,
namely, the A-, B-, C-, and D-approaches. Below we will work with
the B-, C-, D-approaches. We will discuss the relation between geometry and the Hirota's bilinear
method.  Also, we will consider the connection between
supersymmetry, geometry and soliton equations.

\section{Curves and Solitons in 2+1}

In this section, we work with the D-approach. Using this D-approach, we will
establish a connection between curves and (2+1)-dimensional soliton
equations.

\subsection{Some 2-dimensional extensions of the SFE}

According to the D-approach, to establish the connection between (2+1)-dimensional soliton equations
and differential geometry of curves in [1] was  constructed some
two (spatial) dimensional generalizations of the SFE (1a). Here we present
some of them.

\subsubsection{The M-LIX equation}

This equation has the form [1]
$$
\alpha {\bf e}_{1y}=f_{1}{\bf e}_{1x}+
\sum_{j=1}^{n}b_{j}{\bf e}_{1}\wedge \frac{\partial^{j}}{\partial x^{j}}{\bf e}_{1} +
c_{1}{\bf e}_{2}+d_{1}{\bf e}_{3}  \eqno(10a)
$$
$$
\alpha {\bf e}_{2y} =Exercise \quad N1  \eqno(10b)
$$
$$
\alpha {\bf e}_{3y} = Exercise \quad N1 \eqno(10c)
$$
Here the finding of the explicit forms of r.h. of (10b,c) we left
as the exercises (see, the section 7).

\subsubsection{The M-LX equation}

The M-LX equation reads as [1]
$$
\alpha\left ( \begin{array}{ccc}
{\bf e}_{1} \\
{\bf e}_{2} \\
{\bf e}_{3}
\end{array} \right)_{y}= A
\left ( \begin{array}{ccc}
{\bf e}_{1} \\
{\bf e}_{2} \\
{\bf e}_{3}
\end{array} \right)_{x}  + B
\left ( \begin{array}{ccc}
{\bf e}_{1} \\
{\bf e}_{2} \\
{\bf e}_{3}
\end{array} \right)
\eqno(11a)
$$
$$
\left ( \begin{array}{ccc}
{\bf e}_{1} \\
{\bf e}_{2} \\
{\bf e}_{3}
\end{array} \right)_{t}=
\sum_{j=0}^{n}C_{j}\frac{\partial^{j}}{\partial x^{j}}
\left ( \begin{array}{ccc}
{\bf e}_{1} \\
{\bf e}_{2} \\
{\bf e}_{3}
\end{array} \right)
\eqno(11b)
$$
where $A, B, C_{j}$ - some matrices.

\subsubsection{The M-LXI equation}

This extension has the form [1]
$$
\left ( \begin{array}{c}
{\bf e}_{1} \\
{\bf e}_{2} \\
{\bf e}_{3}
\end{array} \right)_{x}= C
\left ( \begin{array}{ccc}
{\bf e}_{1} \\
{\bf e}_{2} \\
{\bf e}_{3}
\end{array} \right), \quad \left ( \begin{array}{ccc}
{\bf e}_{1} \\
{\bf e}_{2} \\
{\bf e}_{3}
\end{array} \right)_{y}= D
\left ( \begin{array}{ccc}
{\bf e}_{1} \\
{\bf e}_{2} \\
{\bf e}_{3}
\end{array} \right)
\eqno(12a)
$$
где
$$
C =
\left ( \begin{array}{ccc}
0             & k     &  0 \\
-\beta k      & 0     & \tau  \\
0             & -\tau & 0
\end{array} \right) ,
\quad
D=
\left ( \begin{array}{ccc}
0            & m_{3}  & -m_{2} \\
-\beta m_{3} & 0      & m_{1} \\
\beta m_{2}  & -m_{1} & 0
\end{array} \right).  \eqno(12b)
$$

\subsubsection{The modified M-LXI equation}

The modified M-LXI (mM-LXI) equation usually we write in the form [1]
$$
\left ( \begin{array}{c}
{\bf e}_{1} \\
{\bf e}_{2} \\
{\bf e}_{3}
\end{array} \right)_{x}= C_{m}
\left ( \begin{array}{ccc}
{\bf e}_{1} \\
{\bf e}_{2} \\
{\bf e}_{3}
\end{array} \right), \quad \left ( \begin{array}{ccc}
{\bf e}_{1} \\
{\bf e}_{2} \\
{\bf e}_{3}
\end{array} \right)_{y}= D_{m}
\left ( \begin{array}{ccc}
{\bf e}_{1} \\
{\bf e}_{2} \\
{\bf e}_{3}
\end{array} \right)
\eqno(13a)
$$
where
$$
C_{m} =
\left ( \begin{array}{ccc}
0             & k     &  -\sigma \\
-\beta k      & 0     & \tau  \\
\beta\sigma             & -\tau & 0
\end{array} \right) ,
\quad
D_{m}= D=
\left ( \begin{array}{ccc}
0            & m_{3}  & -m_{2} \\
-\beta m_{3} & 0      & m_{1} \\
\beta m_{2}  & -m_{1} & 0
\end{array} \right)  \eqno(13b)
$$
and so on [1]. In this paper, we work with the M-LIX, M-LXI and mM-LXI
equations. Note that the M-LXI equation is the particular case of the
mM-LXI eq. as $\sigma=0$.

\subsection{The mM-LXI equation and the mM-LXII equation}

Let us return to the mM-LXI equation (13), which we write in the form
$$
\left ( \begin{array}{c}
{\bf e}_{1} \\
{\bf e}_{2} \\
{\bf e}_{3}
\end{array} \right)_{x}= C_{m}
\left ( \begin{array}{ccc}
{\bf e}_{1} \\
{\bf e}_{2} \\
{\bf e}_{3}
\end{array} \right) \eqno(14a)
$$
$$
\left ( \begin{array}{ccc}
{\bf e}_{1} \\
{\bf e}_{2} \\
{\bf e}_{3}
\end{array} \right)_{y}= D_{m}
\left ( \begin{array}{ccc}
{\bf e}_{1} \\
{\bf e}_{2} \\
{\bf e}_{3}
\end{array} \right)
\eqno(14b)
$$
$$
\left ( \begin{array}{ccc}
{\bf e}_{1} \\
{\bf e}_{2} \\
{\bf e}_{3}
\end{array} \right)_{t}= G
\left ( \begin{array}{ccc}
{\bf e}_{1} \\
{\bf e}_{2} \\
{\bf e}_{3}
\end{array} \right) \eqno(14c)
$$
where
$$
G =
\left ( \begin{array}{ccc}
0       & \omega_{3}  & -\omega_{2} \\
-\beta\omega_{3} & 0      & \omega_{1} \\
\beta\omega_{2}  & -\omega_{1} & 0
\end{array} \right)  \eqno(15)
$$

From (14a,b), we obtain the following mM-LXII equation [1]
$$
C_y - D_x + [C, D] = 0 \eqno (16a)
$$
or
$$
k_{y} - m_{3x} + \sigma m_{1} - \tau m_{2} = 0 \eqno(16b)
$$
$$
\sigma_{y} - m_{2x} + \tau m_{3} - km_{1} =0   \eqno(16c)
$$
$$
\tau_{y} - m_{1x} + \beta (km_{2} - \sigma m_{3}) =0. \eqno(16d)
$$
As $\sigma=0$ the mM-LXII equation reduces to the M-LXII equation [1].
The mM-LXII  equation (16), we can rewrite in form
$$
k_{y} - m_{3x} = \frac{1}{\beta}{\bf e}_{3}\dot ({\bf e}_{3x}\wedge {\bf e}_{3y}) \eqno(17a)
$$
$$
\sigma_{y} - m_{2x} = \frac{1}{\beta}{\bf e}_{2}\dot ({\bf e}_{2x}\wedge {\bf e}_{2y}) \eqno(17b)
$$
$$
\tau_{y} - m_{1x}  = {\bf e}_{1}\dot ({\bf e}_{1x}\wedge {\bf e}_{1y}) \eqno(17c)
$$

Also from (14) we get
$$
k_{t} - \omega_{3x} + \sigma \omega_{1} - \tau \omega_{2} = 0 \eqno (18a)
$$
$$
\sigma_{t} - \omega_{2x} + \tau \omega_{3} - k \omega_{1} = 0 \eqno (18b)
$$
$$
\tau_{t} - \omega_{1x} + \beta (k \omega_{2} - \sigma \omega_{3}) = 0 \eqno (18c)
$$
and
$$
m_{1t} - \omega_{1y} + \beta (m_{3} \omega_{2} - m_{2} \omega_{3}) = 0 \eqno (19a)
$$
$$
m_{2t} - \omega_{2y} + m_{1} \omega_{3} - m_{3} \omega_{1} = 0 \eqno (19b)
$$
$$
m_{3t} - \omega_{3y} + m_{2} \omega_{1} - m_{1} \omega_{2} = 0. \eqno (19c)
$$

\subsection{On the topological invariants}

From  the mM-LXII equation (16) follows
$$
[C, D]_{t} + C_{ty} - D_{tx}  = 0 \eqno (20a)
$$
or
$$
(\sigma m_{1} - \tau m_{2})_{t} + k_{ty} - m_{3tx}  = 0 \eqno(20b)
$$
$$
(\tau m_{3} - km_{1})_{t} + \sigma_{ty} - m_{2tx} = 0   \eqno(20c)
$$
$$
\epsilon (km_{2} - \sigma m_{3})_{t} + \tau_{ty} - m_{1tx} = 0. \eqno(20d)
$$
Hence we get
$$
(\sigma m_1 -\tau m_2)_t +(\sigma \omega_1 -\tau \omega_2)_y -
(m_2 \omega_1 -m_1 \omega_2)_x =0, \eqno (21a)
$$
$$
(\tau m_3 -k m_1)_t +(\tau \omega_3 -k \omega_1)_y -
(m_1 \omega_3 -m_3 \omega_1)_x =0, \eqno (21b)
$$
$$
(k m_2 -\sigma m_3)_t +(k \omega_2 -\sigma \omega_3)_y -
(m_3 \omega_2 -m_2 \omega_3)_x =0. \eqno (21c)
$$
So we have proved the following
\\
{\bf Teorema}: The (2+1)-dimensional nonlinear evolution equations (NLEE)
or dynamical curves which
are given by the mM-LXI equation have the following integrals of motions
$$
K_{1} = \int \int (\kappa m_2 +\sigma m_{3})dxdy, \quad
K_{2} = \int \int (\tau m_2 + \sigma m_{1})dxdy, \quad
K_{3} = \int \int (\tau m_3 -km_{1})dxdy \eqno(22a)
$$
or
$$
K_{1} = \int \int {\bf e}_1({\bf e}_{1x} \wedge {\bf e}_{1y})dxdy \eqno(22b)
$$
$$
K_{2} = \int \int {\bf e}_2({\bf e}_{2x} \wedge {\bf e}_{2y})dxdy \eqno(22c)
$$
$$
K_{3} = \int \int {\bf e}_3({\bf e}_{3x} \wedge {\bf e}_{3y})dxdy. \eqno(22d)
$$

So we have the following three topological invariants
$$
Q_{1} =\frac{1}{4\pi} \int \int {\bf e}_{1}\dot ({\bf e}_{1x}\wedge {\bf e}_{1y})dxdy \eqno(23a)
$$
$$
Q_{2} =\frac{1}{4\pi} \int \int {\bf e}_{2}\dot ({\bf e}_{2x}\wedge {\bf e}_{2y})dxdy \eqno(23b)
$$
$$
Q_{3} =\frac{1}{4\pi} \int \int {\bf e}_{3}\dot ({\bf e}_{3x}\wedge {\bf e}_{3y})dxdy \eqno(23c)
$$
We note that may be not all of these topological invariants are independent.

\subsection{The M-LXI equation and Soliton equations in 2+1}

In this section we will establish the connection between the
M-LXI equation (12) and soliton equations in 2+1 dimensions.
Let us, we assume
$$
{\bf e}_{1} \equiv {\bf S} \eqno(24)
$$
Moreover we introduce two complex functions
$q, p$ according to the following expressions
$$
q = a_{1}e^{ib_{1}}, \quad p=a_{2}e^{ib_{2}}  \eqno(25)
$$
where $a_{j}, b_{j}$ are real functions. Now we ready to consider some
examples.

\subsubsection{The Ishimori equation}

The Ishimori equation (IE) reads as [7]
$$
{\bf S}_{t}  =  {\bf S}\wedge ({\bf S}_{xx} +\alpha^2 {\bf S}_{yy})+
u_x{\bf S}_{y}+u_y{\bf S}_{x} \eqno (26a)
$$
$$
u_{xx}-\alpha^2 u_{yy}  = -2\alpha^2 {\bf S}\cdot ({\bf S}_{x}\wedge
                        {\bf S}_{y}).   \eqno (26b)
$$
In this case we have
$$
m_{1}=\partial_{x}^{-1}[\tau_{y}-\frac{\epsilon}{2\alpha^2}M_2^{Ish}u] \eqno(27a)
$$
$$
m_{2}= -\frac{1}{2\alpha^2 k}M_2^{Ish}u \eqno (27b)
$$
$$
m_{3} =\partial_{x}^{-1}[k_y +\frac{\tau}{2\alpha^2 k}M_2^{Ish}u] \eqno(27c)
$$
and
$$
\omega_{1} = \frac{1}{k}[-\omega_{2x}+\tau\omega_{3}]
                 \eqno (28a)
$$
$$
\omega_{2}= -k_{x}-
 \alpha^{2}(m_{3y}+m_{2}m_{1})+im_{2}u_{x}
\eqno (28b)
$$
$$
\omega_{3}= -k \tau+\alpha^{2}(m_{2y}-m_{3}m_{1})
+ik u_{y}+im_{3}u_{x}.
\eqno (28c)
$$
$$
M_2^{Ish}=M_2|_{a=b=-\frac{1}{2}}.
$$
Functions $q, p$ are given by (25) with
$$
a_{1}^2 =a_{1}^{\prime^{2}}=\frac{1}{4}k^2+
\frac{|\alpha|^2}{4}(m_3^2 +m_2^2)-\frac{1}{2}\alpha_{R}km_3-
\frac{1}{2}\alpha_{I}km_2
  \eqno(29a)
$$
$$
b_1 =\partial_{x}^{-1}\{-\frac{\gamma_1}{2ia_1^{\prime^{2}}}-(\bar A-A+D-\bar D)\}  \eqno(29b)
$$
$$
a_2^2=a_{2}^{\prime^{2}}=\frac{1}{4}k^2+
\frac{|\alpha|^2}{4}(m_3^2 +m_2^2)+\frac{1}{2}\alpha_{R}km_3
-\frac{1}{2}\alpha_{I}km_2
  \eqno(29c)
$$
$$
b_{2} =\partial_{x}^{-1}\{-\frac{\gamma_2}{2ia_2^{\prime^{2}}}-(A-\bar A+\bar D-D)\}  \eqno(29d)
$$
where
$$
\gamma_1=i\{\frac{1}{2}k^{2}\tau+
\frac{|\alpha|^2}{2}(m_3km_1+m_2k_y)-
$$
$$
\frac{1}{2}\alpha_{R}(k^{2}m_1+m_3k\tau+
m_2k_x)
+\frac{1}{2}\alpha_{I}[k(2k_y-m_{3x})-
k_x m_3]\}. \eqno(30a)
$$
$$
\gamma_2=-i\{\frac{1}{2}k^{2}\tau+
\frac{|\alpha|^2}{2}(m_3km_1+m_2 k_y)+
$$
$$
\frac{1}{2}\alpha_{R}(k^{2}m_1+m_3k\tau+
m_2k_x )
+\frac{1}{2}\alpha_{I}[k(2k_y-m_{3x})-
k_x m_3]\}. \eqno(30b)
$$
Here $\alpha=\alpha_{R}+i\alpha_{I}$. In this case, $q,p$ satisfy the following
DS  equation
$$
iq_t  + q_{xx}+\alpha^{2}q_{yy} + vq = 0 \eqno (31a)
$$
$$
-ip_t +  p_{xx}+\alpha^{2}p_{yy} + vp = 0 \eqno (31b)
$$
$$
v_{xx}-\alpha^{2}v_{yy} + 2[(p q)_{xx}+\alpha^{2}(p q)_{yy}] = 0.
\eqno (31c)
$$
So we have proved that theIE (26) and the (31) are L-equivalent to each other.
As well known that these equations are G-equivalent to each other [5].
Note that the IE contains two reductions:
the Ishimori I equation as $\alpha_{R}=1, \alpha_{I}=0$ and
the Ishimori II equation as $\alpha_{R}=0, \alpha_{I}=1$.
The corresponding versions of the DS equation (31), we obtain as
the corresponding values of the parameter $\alpha$ [1].

\subsubsection{The Myrzakulov IX  equation}

Now we find the connection between the Myrzakulov IX (M-IX) equation
and the curves (the M-LXI equation). The M-IX equation reads as
$$
{\bf S}_t = {\bf S} \wedge M_1{\bf S}+A_2{\bf S}_x+A_1{\bf S}_y  \eqno(32a)
$$
$$
M_2u=2\alpha^{2} {\bf S}({\bf S}_x \wedge {\bf S}_y) \eqno(32b)
$$
where $ \alpha,b,a  $=  consts and
$$
M_1= \alpha ^2\frac{\partial ^2}{\partial y^2}+4\alpha (b-a)\frac{\partial^2}
   {\partial x \partial y}+4(a^2-2ab-b)\frac{\partial^2}{\partial x^2},
$$
$$
M_2=\alpha^2\frac{\partial^2}{\partial y^2} -2\alpha(2a+1)\frac{\partial^2}
   {\partial x \partial y}+4a(a+1)\frac{\partial^2}{\partial x^2},
$$
$$
A_1=i\{\alpha (2b+1)u_y - 2(2ab+a+b)u_{x}\},
$$
$$
A_2=i\{4\alpha^{-1}(2a^2b+a^2+2ab+b)u_x - 2(2ab+a+b)u_{y}\}.
$$
The M-IX equation was introduced in [1] and is integrable. It admits several
integrable reductions:
\\
1) the Ishimori equation as $a=b=-\frac{1}{2}$
\\
2) the M-VIII equation as $a=b=-1$
\\
and so on [1].
In this case we have
$$
m_{1}=\partial_{x}^{-1}[\tau_{y}-\frac{\beta}{2\alpha^2}M_2 u] \eqno(33a)
$$
$$
m_{2}=-\frac{1}{2\alpha^2 k}M_2 u  \eqno (33b)
$$
$$
m_{3}=\partial_{x}^{-1}[k_y +\frac{\tau}{2\alpha^2 k}M_2 u]  \eqno(33c)
$$
and
$$
\omega_{1} = \frac{1}{k}[-\omega_{2x}+\tau\omega_{3}],
                 \eqno (34a)
$$
$$
\omega_{2}= -4(a^{2}-2ab-b)k_{x}-
4\alpha (b-a)k_{y} -\alpha^{2}(m_{3y}+m_{2}m_{1})+m_{2}A_{1}
\eqno (34b)
$$
$$
\omega_{3}= -4(a^{2}-2ab-b)k \tau-
4\alpha (b-a)k m_{1}+\alpha^{2}(m_{2y}-m_{3}m_{1})
+k A_{2}+m_{3}A_{1}
\eqno (34c)
$$
Functions $q, p$ are given by (25) with
$$
a_{1}^2 =\frac{|a|^2}{|b|^2}a_{1}^{\prime^{2}}=\frac{|a|^2}{|b|^2}\{(l+1)^2k^2
+\frac{|\alpha|^2}{4}(m_3^2 +m_2^2)-(l+1)\alpha_{R}km_3-
(l+1)\alpha_{I}km_2\}
  \eqno(35a)
$$
$$
b_{1} =\partial_{x}^{-1}\{-\frac{\gamma_1}{2ia_1^{\prime^{2}}}-(\bar A-A+D-\bar D)\}  \eqno(35b)
$$
$$
a_{2}^2 =\frac{|b|^2}{|a|^2}a_{2}^{\prime^{2}}=\frac{|b|^2}{|a|^2}\{l^2k^2
+\frac{|\alpha|^2}{4}(m_3^2 +m_2^2)-l\alpha_{R}km_3+
l\alpha_{I}km_2\}
  \eqno(35c)
$$
$$
b_{2} =\partial_{x}^{-1}\{-\frac{\gamma_2}{2ia_2^{\prime^{2}}}-(A-\bar A+\bar D-D)  \eqno(2.9a)
$$
where
$$
\gamma_1=i\{2(l+1)^2k^{2}\tau+\frac{|\alpha|^2}{2}(m_3km_1+m_2k_y)-
$$
$$
(l+1)\alpha_{R}[k^{2}m_1+m_3k\tau+
m_2k_x]+(l+1)\alpha_{I}[k(2k_y-m_{3x})-
k_x m_3]\} \eqno(36a)
$$
$$
\gamma_2=-i\{2l^2k^{2}\tau+
\frac{|\alpha|^2}{2}(m_3km_1+m_2k_y)-
$$
$$
l\alpha_{R}(k^{2}m_1+m_3k\tau+
m_2k_x)-l\alpha_{I}[k(2k_y-m_{3x})-
k_x m_3]\}. \eqno(36b)
$$
Here $\alpha=\alpha_{R}+i\alpha_{I}$. In this case, $q,p$
satisfy the following Zakharov equation [4]
$$
iq_t+M_{1}q+vq=0 \eqno(37a)
$$
$$
ip_t-M_{1}p-vp=0 \eqno(37b)
$$
$$
M_{2}v=-2M_{1}(pq) \eqno(37c)
$$

As well known the M-IX equation admits several reductions:
1)  the M-IXA equation as $\alpha_{R}=1, \alpha_{I}=0$;
2)  the M-IXB equation as $\alpha_{R}=0, \alpha_{I}=1$;
3) the M-VIII equation as $a=b=1$
4) the IE $a=b=-\frac{1}{2}$
and so on. The corresponding versions of the ZE (9), we obtain as
the corresponding values of the parameter $\alpha$.

\subsection{The modified M-LXI equation and Soliton equations in 2+1}

In this section we will establish the connection between the modified
M-LXI equation (14) and soliton equations in 2+1 dimensions. As above we assume
$$
{\bf e}_{1} \equiv {\bf S} \eqno(38)
$$
and
$$
q = a_{1}e^{ib_{1}}, \quad p=a_{2}e^{ib_{2}}  \eqno(39)
$$
where $a_{j}, b_{j}$ are as and above, real functions. Examples.

\subsubsection{The Ishimori equation}

Consider the IE (26).
For this equation we obtain
$$
m_{1}=\partial_{x}^{-1}[\tau_{y}-\frac{\beta}{2\alpha^2}M_2^{Ish}u] \eqno(40a)
$$
$$
m_{2}=\frac{\sigma}{k}m_3 -\frac{1}{2\alpha^2 k}M_2^{Ish}u  \eqno (40b)
$$
$$
m_{3x} +\frac{\tau\sigma}{k}m_3=k_y +\sigma\partial_x^{-1}[\tau_y -
\frac{\epsilon}{2\alpha^2}M_2^{Ish}u]+\frac{\tau}{2\alpha^2 k}M_2^{Ish}u \eqno(40c)
$$
and
$$
\omega_{1} = \frac{1}{k}[\sigma_{t}-\omega_{2x}+\tau\omega_{3}]
                 \eqno (41a)
$$
$$
\omega_{2}= -(k_{x}+\sigma \tau)-
\alpha^{2}(m_{3y}+m_{2}m_{1})+i\sigma u_{y}+im_{2}u_{x}
\eqno (41b)
$$
$$
\omega_{3}= (\sigma_{x}-k \tau)+
\alpha^{2}(m_{2y}-m_{3}m_{1})
+ik u_{y}+im_{3}u_{x}.
\eqno (41c)
$$
Functions $q, p$ are given by (39) with
$$
a_{1}^2 =a_{1}^{\prime^{2}}=\frac{1}{4}(k+
\sigma^2)+\frac{|\alpha|^2}{4}(m_3^2 +m_2^2)-frac{1}{2}\alpha_{R}(km_3+\sigma m_2)-
\frac{1}{2}\alpha_{I}(km_2+\sigma m_3)
  \eqno(42a)
$$
$$
b_{1} =\partial_{x}^{-1}\{-\frac{\gamma_1}{2ia_1^{\prime^{2}}}-(\bar A-A+D-\bar D)\}  \eqno(42b)
$$
$$
a_2^2=a_{2}^{\prime^{2}}=\frac{1}{4}(k^2+
\sigma^2)+\frac{|\alpha|^2}{4}(m_3^2 +m_2^2)+\frac{1}{2}\alpha_{R}(km_3+\sigma m_2)
-\frac{1}{2}\alpha_{I}(km_2+\sigma m_3)\}
  \eqno(42c)
$$
$$
b_{2} =\partial_{x}^{-1}\{-\frac{\gamma_2}{2ia_2^{\prime^{2}}}-(A-\bar A+\bar D-D)\}  \eqno(42d)
$$
where
$$
\gamma_1=i\{\frac{1}{2}[k(k\tau-\sigma_x)+\sigma(\sigma\tau+k_x)]+
\frac{|\alpha|^2}{2}[m_3(km_1-\sigma_y)+m_2(\sigma m_1 +k_y)]-
$$
$$
\frac{1}{2}\alpha_{R}[k(km_1-\sigma_y)+\sigma(\sigma m_1+k_y)+m_3(k\tau-\sigma_x)+
m_2(\sigma\tau+k_x)]+
$$
$$
\frac{1}{2}\alpha_{I}[k(2k_y-m_{3x})+\sigma(2\sigma_y-m_{2x})-
k_x m_3-\sigma_x m_2]\} \eqno(43a)
$$
$$
\gamma_2=-i\{\frac{1}{2}[k(k\tau-\sigma_x)+\sigma(\sigma\tau+k_x)]+
\frac{|\alpha|^2}{2}[m_3(km_1-\sigma_y)+m_2(\sigma m_1 +k_y)]+
$$
$$
\frac{1}{2}\alpha_{R}[k(km_1-\sigma_y)+\sigma(\sigma m_1+k_y)+m_3(k\tau-\sigma_x)+
m_2(\sigma\tau+k_x)]+
$$
$$
\frac{1}{2}\alpha_{I}[k(2k_y-m_{3x})+\sigma(2\sigma_y-m_{2x})-
k_x m_3-\sigma_x m_2]\}. \eqno(43b)
$$
Here
$$
\alpha=\alpha_{R}+i\alpha_{I}, \quad A=\frac{i}{4}[u_{y}-
\frac{2a}{\alpha}u_{x}],
\quad D=\frac{i}{4}[\frac{(2a+1)}{\alpha}u_{x}-u_{y}].
$$
 In this case, $q,p$ satisfy the
DS  equation (31).

The Ishimori I and DS I equations, we get as $\alpha_R=1, \alpha_I=0.$
The Ishimori II and DS II equations we obtain from these results as
$\alpha_R=0, \alpha_I=1.$ Details, you can find in [1].

\subsubsection{The Myrzakulov IX  equation}

Now let us establish the connection between the M-IX equation (32)
and the mM-LXI equation (14). From (32) and (14) we get
$$
m_{1}=\partial_{x}^{-1}[\tau_{y}-\frac{\beta}{2\alpha^2}M_2 u] \eqno(44a)
$$
$$
m_{2}=\frac{\sigma}{k}m_3 -\frac{1}{2\alpha^2 k}M_2 u  \eqno (44b)
$$
$$
m_{3x} +\frac{\tau\sigma}{k}m_3=k_y +\sigma\partial_x^{-1}[\tau_y -
\frac{\epsilon}{2\alpha^2}M_2 u]+\frac{\tau}{2\alpha^2 k}M_2 u  \eqno(44c)
$$
and
$$
\omega_{1} = \frac{1}{k}[\sigma_{t}-\omega_{2x}+\tau\omega_{3}],
                 \eqno (45a)
$$
$$
\omega_{2}= -4(a^{2}-2ab-b)(k_{x}+\sigma \tau)-
4\alpha (b-a)(k_{y}+\sigma m_{1}) -\alpha^{2}(m_{3y}+m_{2}m_{1})+\sigma A_{2}+m_{2}A_{1}
\eqno (45b)
$$
$$
\omega_{3}= 4(a^{2}-2ab-b)(\sigma_{x}-k \tau)+
4\alpha (b-a)(\sigma_{y}-k m_{1}) +\alpha^{2}(m_{2y}-m_{3}m_{1})
+k A_{2}+m_{3}A_{1}
\eqno (45c)
$$
Functions $q, p$ are given by (39) with
$$
a_{1}^2 =\frac{|a|^2}{|b|^2}a_{1}^{\prime^{2}}=\frac{|a|^2}{|b|^2}\{(l+1)^2(k+
\sigma^2)+\frac{|\alpha|^2}{4}(m_3^2 +m_2^2)-(l+1)\alpha_{R}(km_3+\sigma m_2)-
(l+1)\alpha_{I}(km_2+\sigma m_3)\}
  \eqno(46a)
$$
$$
b_{1} =\partial_{x}^{-1}\{-\frac{\gamma_1}{2ia_1^{\prime^{2}}}-(\bar A-A+D-\bar D)\}  \eqno(46b)
$$
$$
a_{2}^2 =\frac{|b|^2}{|a|^2}a_{2}^{\prime^{2}}=\frac{|b|^2}{|a|^2}\{l^2(k^2+
\sigma^2)+\frac{|\alpha|^2}{4}(m_3^2 +m_2^2)-l\alpha_{R}(km_3+\sigma m_2)+
l\alpha_{I}(km_2+\sigma m_3)\}
  \eqno(46c)
$$
$$
b_{2} =\partial_{x}^{-1}\{-\frac{\gamma_2}{2ia_2^{\prime^{2}}}-(A-\bar A+\bar D-D)\}  \eqno(46d)
$$
where
$$
\gamma_1=i\{2(l+1)^2[k(k\tau-\sigma_x)+\sigma(\sigma\tau+k_x)]+
\frac{|\alpha|^2}{2}[m_3(km_1-\sigma_y)+m_2(\sigma m_1 +k_y)]-
$$
$$
(l+1)\alpha_{R}[k(km_1-\sigma_y)+\sigma(\sigma m_1+k_y)+m_3(k\tau-\sigma_x)+
m_2(\sigma\tau+k_x)]+
$$
$$
(l+1)\alpha_{I}[k(2k_y-m_{3x})+\sigma(2\sigma_y-m_{2x})-
k_x m_3-\sigma_x m_2]\} \eqno (47a)
$$
$$
\gamma_2=-i\{2l^2[k(k\tau-\sigma_x)+\sigma(\sigma\tau+k_x)]+
\frac{|\alpha|^2}{2}[m_3(km_1-\sigma_y)+m_2(\sigma m_1 +k_y)]-
$$
$$
l\alpha_{R}[k(km_1-\sigma_y)+\sigma(\sigma m_1+k_y)+m_3(k\tau-\sigma_x)+
m_2(\sigma\tau+k_x)]-
$$
$$
l\alpha_{I}[k(2k_y-m_{3x})+\sigma(2\sigma_y-m_{2x})-
k_x m_3-\sigma_x m_2]\}. \eqno (47b)
$$
Directly calculation show that $q,p$ satisfy the ZE (37).

These results gives: 1) as $\alpha_R=1, \alpha_I=0$ the M-IXA equation;
2) as $\alpha_R=0, \alpha_I=1$ the M-IXB equation; 3) as $a=b=-\frac{1}{2},
\alpha_R=1, \alpha_I=0$ the Ishimori I and DS I equations; 4) as $a=b=-\frac{1}{2},
\alpha_R=0, \alpha_I=1$ the Ishimori II and DS II equations; 5) as $a=b-1$
the M-VIII and corresponding Zakharov equations; and so on [1].

\subsection{The M-LIX equation and Soliton equations in 2+1}

Now let us consider the connection between the M-LIX equation and
(2+1)-dimensional soliton equations. Mention that the M-LIX equation is
one of (2+1)-dimensional extensions of the SFE (1a).
As example, let us consider the connection between the M-LIX equation
and the M-IX
equation (32). Let the M-LIX equation has the form [1]
$$
\alpha {\bf e}_{1y}=\frac{2a+1}{2}{\bf e}_{1x}+
\frac{i}{2}{\bf e}_{1}\wedge{\bf e}_{1x} +
+i(q+p){\bf e}_{2}+(q-p){\bf e}_{3}  \eqno(48a)
$$
$$
\alpha {\bf e}_{2y} = Exercise \quad N1 \eqno(48b)
$$
$$
\alpha {\bf e}_{3y} = Exercise \quad N1. \eqno(48c)
$$
In terms of matrix this equation we can write in the form
$$
\alpha \hat  e_{1y} = \frac{2a+1}{2}\hat e_{1x}
+ \frac{1}{4}[\hat e_{1},\hat e_{1x}] +
i(q+p)\hat e_{2}+(q-p)\hat e_{3}  \eqno(49a)
$$
$$
\alpha \hat  e_{2y} =  Exercise \quad N1  \eqno(49b)
$$
$$
\alpha \hat  e_{3y} =  Exercise \quad N1  \eqno(49c)
$$
where
$$
\hat  e_{1} = g^{-1}\sigma_{3}g, \quad \hat e_{2}=g^{-1}\sigma_{2}g,
\quad \hat e_{3} = g^{-1}\sigma_{1}g  \eqno(50)
$$
Here  $\sigma_{j}$ are Pauli matrices
$$
\sigma_{1} =
\left ( \begin{array}{cc}
0       &  1 \\
1       &  0
\end{array} \right) , \quad
\sigma_{2} =
\left ( \begin{array}{cc}
0       &  -i \\
i       &  0
\end{array} \right) , \quad
\sigma_{3} =
\left ( \begin{array}{cc}
1       &  0 \\
0       &  -1
\end{array} \right)   \eqno(51)
$$
So we have
$$
\sigma_{1}\sigma_{2} = i\sigma_{3} = -\sigma_{2}\sigma_{1}, \quad
\sigma_{1}\sigma_{3} = -i\sigma_{2} = -\sigma_{3}\sigma_{1}, \quad
\sigma_{3}\sigma_{2} = -i\sigma_{1} = -\sigma_{2}\sigma_{3}
\eqno(52a)
$$
and
$$
\sigma_{j}^{2} =I = diag(1,1). \eqno(52b)
$$
Equations (49) we can rewrite in the form
$$
[\sigma_{3},B_{0}] = i(q+p)\sigma_{2}+(q-p)\sigma_{1}   \eqno(53a)
$$
$$
[\sigma_{2},B_{0}] = -i(q+p)\sigma_{3}   \eqno(53b)
$$
$$
[\sigma_{1},B_{0}] =-(q-p)\sigma_{3}   \eqno(53c)
$$
where
$$
B_{0} = \alpha g_{y}g^{-1} - B_{1}g_{x}g^{-1}, \quad
B_{1}=\frac{2a+1}{2}I+\frac{1}{2}\sigma_{3}    \eqno(54)
$$
Hence we get
$$
B_{0} =
\left ( \begin{array}{cc}
0       &  q \\
p       &  0
\end{array} \right).   \eqno(55)
$$
Thus the matrix-function  $g$ satisfies the equations
$$
\alpha g_{y}=  B_{1}g_{x} + B_{0}g. \eqno(56)
$$

To find the time evolution of matrices  $\hat e_{j}$ or vectors
${\bf e}_{j}$, we require that the matrix  $\hat e_{1}$
saisfy the  M-IX equation, i.e.
$$
i\hat  e_{1t} = \frac{1}{2}[\hat e_{1}, M_{1}\hat e_{1}] +A_{1}\hat e_{1y}+A_{2}\hat e_{1x}
\eqno(57a)
$$
$$
M_{2}u = \frac{\alpha^{2}}{2i}tr(\hat e_{1}([\hat e_{1x},\hat e_{1y}])
\eqno(57b)
$$

From these informations we find the time evolution of matrices
 $\hat e_{2}, \hat e_{3}$. So after some algebra we obtain
$$
[\sigma_{3},C_{0}] = i(c_{12}+c_{21})\sigma_{2}+(c_{12}-c_{21})\sigma_{1}   \eqno(58a)
$$
$$
[\sigma_{2},C_{0}] = i(c_{11}-c_{22})\sigma_{1}-i(c_{12}+c_{21})\sigma_{3}   \eqno(58b)
$$
$$
[\sigma_{1},C_{0}] =-i(c_{11}-c_{22})\sigma_{2}-(c_{1@}-c_{21})\sigma_{3}   \eqno(58c)
$$
where
$$
C_{0} = g_{t}g^{-1} - 2iC_{2}g_{xx}g^{-1}-C_{1}g_{x}g^{-1}, \quad
C_{2}=\frac{2b+1}{2}I+\frac{1}{2}\sigma_{3},
\quad C_{1}=iB_{0}.    \eqno(59)
$$
Hence we get
$$
C_{0} =
\left ( \begin{array}{cc}
c_{11}      &  c_{12} \\
c_{21}       &  c_{22}
\end{array} \right)   \eqno(60)
$$
with
$$
c_{12}=i(2b-a+1)q_{x}+i\alpha q_{y},\quad c_{21}=i(a-2b)q_{x}-i\alpha p_{y} \eqno(61a)
$$
and $c_{jj}$ are the solutions of the following equations
$$
(a+1)c_{11x}-\alpha c_{11y}=i[(2b-a+1)(pq)_{x}+\alpha(pq)_{y}]  \eqno(61b)
$$
$$
ac_{22x}-\alpha c_{22y}=i[(a-2b)(pq)_{x}-\alpha (pq)_{y}].       \eqno(61c)
$$
So that the matrix $g$ satisfies the equation
$$
g_{t}=  2C_{2}g_{xx} + C_{1}g_{x} +C_{0}g. \eqno(62)
$$

So we have identified the curve,
given by the M-LIX equation (48) with
the M-IX equation (32). On the other hand, the compatibilty
condition of equations (56) and (62) is equivalent to the ZE (37).
So that we have also established the connection between the  curve
(the M-LIX equation) and the ZE. And we have shown, once more that
the M-IX equation (32) and the ZE (37) are L-equivalent to each other. Finally
we note as $a=b=-\frac{1}{2}$ from these results follows the corresponding
connection between the M-LIX, Ishimori and DS equations [1]. And
as $a=b=-1$ we get the relation between the
M-VIII, M-LIX and other Zakharov equations (for details, see [1]).

\subsection{Spin systems as reductions of the M-0 equation}

Consider the (2+1)-dimensional M-0 equation [1]
$$
{\bf S}_{t} = a_{12} {\bf e}_{2} + a_{13}{\bf e}_{3}, \quad
{\bf S}_{x} = b_{12} {\bf e}_{2} + b_{13}{\bf e}_{3}, \quad
{\bf S}_{y} = c_{12} {\bf e}_{2} + c_{13}{\bf e}_{3} \eqno(63)
$$
where
$$
{\bf e}_{2} = \frac{c_{13}}{\triangle}{\bf S}_{x} -
\frac{b_{13}}{\triangle}{\bf S}_{y}, \quad
{\bf e}_{3} = -\frac{c_{12}}{\triangle}{\bf S}_{x} +
\frac{b_{12}}{\triangle}{\bf S}_{y}, \quad \triangle =
b_{12}c_{13}- b_{13}c_{12}. \eqno(64)
$$
All known spin systems (integrable and nonintegrable) in 2+1 dimensions
are the particular reductions of the M-0 equation (63). In particular,
the IE (26) is the integrable reduction of equation (63). In this case,
we have
$$
a_{12} = \omega_{3},  a_{13}= -\omega_{2}, b_{12}= k, b_{13}=  -\sigma,
 c_{12}= m_{3},  c_{13}= -m_{2}.  \eqno(65)
$$
Sometimes we use the following form of the M-0 equation [1]
$$
{\bf S}_{t} = d_{2} {\bf S}_{x} + d_{3}{\bf S}_{y} \eqno(66)
$$
with
$$
d_{2}= \frac{a_{12}c_{13}-a_{13}c_{12}}{\triangle}, \quad
d_{3}= \frac{a_{12}b_{13}-a_{13}b_{12}}{\triangle}. \eqno(67)
$$

\section{Surfaces and Solitons in 2+1}

\subsection{The M-LVIII equation and Soliton equationa in 2+1}

In the C-approach [1], our starting point is  the following (2+1)-dimensional
M-LVIII equation [1]
$$
{\bf r}_{t} = \Upsilon_{1} {\bf r}_{x} + \Upsilon_{2} {\bf r}_{y}
+ \Upsilon_{3}{\bf n}   \eqno(68a)
$$
$$
{\bf r}_{xx} = \Gamma^{1}_{11} {\bf r}_{x} + \Gamma^{2}_{11} {\bf r}_{y}
+ L{\bf n}   \eqno(68b)
$$
$$
{\bf r}_{xy} = \Gamma^{1}_{12} {\bf r}_{x} + \Gamma^{2}_{12} {\bf r}_{y}
+ M{\bf n}   \eqno(68c)
$$
$$
{\bf r}_{yy} = \Gamma^{1}_{22} {\bf r}_{x} + \Gamma^{2}_{22} {\bf r}_{y}
+ N{\bf n}   \eqno(68d)
$$
$$
{\bf n}_{x} = p_{11} {\bf r}_{x} + p_{12} {\bf r}_{y}   \eqno(68e)
$$
$$
{\bf n}_{y} = p_{21} {\bf r}_{x} + p_{22} {\bf r}_{y}.   \eqno(68f)
$$
This equation admits several integrable reductions. Practically, all
integrable spin systems in 2+1 dimensions are some integrable reductions
of the M-LVIII equation (68).

\subsection{The M-LXIII equation and Soliton equationa in 2+1}

Sometimes it is convenient to work using the B-approach. In this approach
the starting equation is the  following M-LXIII equation [1]
$$
{\bf r}_{tx} = \Gamma^{1}_{01} {\bf r}_{x} + \Gamma^{2}_{01} {\bf r}_{y}
+ \Gamma^{3}_{01}{\bf n}   \eqno(69a)
$$
$$
{\bf r}_{ty} = \Gamma^{1}_{02} {\bf r}_{x} + \Gamma^{2}_{02} {\bf r}_{y}
+ \Gamma^{3}_{02}{\bf n}   \eqno(69b)
$$
$$
{\bf r}_{xx} = \Gamma^{1}_{11} {\bf r}_{x} + \Gamma^{2}_{11} {\bf r}_{y}
+ L{\bf n}   \eqno(69c)
$$
$$
{\bf r}_{xy} = \Gamma^{1}_{12} {\bf r}_{x} + \Gamma^{2}_{12} {\bf r}_{y}
+ M{\bf n}   \eqno(69d)
$$
$$
{\bf r}_{yy} = \Gamma^{1}_{22} {\bf r}_{x} + \Gamma^{2}_{22} {\bf r}_{y}
+ N{\bf n}   \eqno(69e)
$$
$$
{\bf n}_{t} = p_{01} {\bf r}_{x} + p_{02} {\bf r}_{y}   \eqno(69f)
$$
$$
{\bf n}_{x} = p_{11} {\bf r}_{x} + p_{12} {\bf r}_{y}   \eqno(69g)
$$
$$
{\bf n}_{y} = p_{21} {\bf r}_{x} + p_{22} {\bf r}_{y}.   \eqno(69h)
$$
This equation follows from the M-LVIII equation (68) under the following conditions
$$
\Gamma^{1}_{01}= \Upsilon_{1x}+ \Upsilon_{1} \Gamma^{1}_{11}
+ \Upsilon_{2}\Gamma^{1}_{12}+\Upsilon_{3}p_{11}
$$
$$
\Gamma^{1}_{02}= \Upsilon_{2x}+ \Upsilon_{1} \Gamma^{2}_{11}
+ \Upsilon_{2}\Gamma^{2}_{12}+\Upsilon_{3}p_{12}
$$
$$
\Gamma^{1}_{03}= \Upsilon_{3x}+ \Upsilon_{1} L
+ \Upsilon_{2}M
$$
$$
p_{01}= \frac{F \Gamma^{3}_{02}}{\Lambda}, \quad
p_{02}= -\frac{E \Gamma^{3}_{02}}{\Lambda}, \quad
\Lambda = EG-F^{2} \eqno(70)
$$
Note that the M-LXIII equation (69) usually we use in the following form
$$
Z_{x} = AZ  \eqno(71a)
$$
$$
Z_{y} = BZ  \eqno(71b)
$$
$$
Z_{t} = CZ  \eqno(71c)
$$
where
$Z=({\bf r}_{x}, {\bf r}_{y}, {\bf n})^{t}$ and
$$
A =
\left ( \begin{array}{ccc}
\Gamma^{1}_{11} & \Gamma^{2}_{11} & L \\
\Gamma^{1}_{12} & \Gamma^{2}_{12} & M \\
p_{11}           & p_{12}           & 0
\end{array} \right) ,  \quad
B =
\left ( \begin{array}{ccc}
\Gamma^{1}_{12} & \Gamma^{2}_{12} & M \\
\Gamma^{1}_{22} & \Gamma^{2}_{22} & N \\
p_{21}          & p_{22}          & 0
\end{array} \right), \quad
C =
\left ( \begin{array}{ccc}
\Gamma^{1}_{01} & \Gamma^{2}_{01} & \Gamma^{3}_{01} \\
\Gamma^{1}_{02} & \Gamma^{2}_{02} & \Gamma^{3}_{02} \\
\Gamma^{1}_{03} & \Gamma^{2}_{03} & 0
\end{array} \right) . \eqno(72)
$$

\subsection{The M-LXIV equation}

In this subsection we derive  the M-LXIV equation, which express some
 relations between
coefficients of the M-LXIII equation (69) or (71). From (71) we have
$$
A_{y}-B_{x}+[A,B]=0 \eqno(73a)
$$
$$
A_{t}-C_{x}+[A,C]=0 \eqno(73b)
$$
$$
B_{t}-C_{y}+[B,C]=0 \eqno(73c)
$$
It is the M-LXIV equation. These equations are equivalent the relations
$$
{\bf r}_{yxx} = {\bf r}_{xxy}, \quad
{\bf r}_{yyx}={\bf r}_{xyy}  \eqno(74a)
$$
$$
{\bf r}_{txx} = {\bf r}_{xxt}, \quad
{\bf r}_{txy}={\bf r}_{xyt}, \quad {\bf r}_{tyy}={\bf r}_{yyt}.  \eqno(74b)
$$
Note that (73a) is the well known Codazzi-Mainardi-Peterson equation (CMPE).

\subsection{Orthogonal basis and LR of the M-LXIV equation}

Let us introduce the orthogonal trihedral
$$
{\bf e}_{1} = \frac{{\bf r}_{x}}{\sqrt E}, \quad
{\bf e}_{2} = {\bf n},
\quad {\bf  e}_{3} = {\bf e}_{1} \wedge {\bf  e}_{2}.    \eqno(75)
$$

Let ${\bf e}_{1}^{2}=\beta = \pm 1, {\bf e}_{2}^{2}={\bf e}_{3}^{2}= 1$.
Then these vectors satisfy the following  equations
$$
\left ( \begin{array}{ccc}
{\bf e}_{1} \\
{\bf e}_{2} \\
{\bf e}_{3}
\end{array} \right)_{x}=\frac{1}{\sqrt{E}}
\left ( \begin{array}{ccc}
0             & L   & -\frac{\Lambda}{\sqrt{E}}\Gamma^{2}_{11} \\
-\beta L      & 0     & -\Lambda p_{12}  \\
\frac{\beta\Lambda}{\sqrt{E}}\Gamma^{2}_{11}& \Lambda  p_{12} & 0
\end{array} \right)
\left ( \begin{array}{c}
{\bf e}_{1} \\
{\bf e}_{2} \\
{\bf e}_{3}
\end{array} \right)
\eqno(76a)
$$
$$
\left ( \begin{array}{ccc}
{\bf e}_{1} \\
{\bf e}_{2} \\
{\bf e}_{3}
\end{array} \right)_{y}=\frac{1}{\sqrt{E}}
\left ( \begin{array}{ccc}
0             & M   & -\frac{\Lambda}{\sqrt{E}}\Gamma^{2}_{12} \\
-\beta M      & 0     & -\Lambda p_{22}  \\
\frac{\beta\Lambda}{\sqrt{E}}\Gamma^{2}_{12}& \Lambda  p_{22} & 0
\end{array} \right)
\left ( \begin{array}{c}
{\bf e}_{1} \\
{\bf e}_{2} \\
{\bf e}_{3}
\end{array} \right)
\eqno(76b)
$$
$$
\left ( \begin{array}{ccc}
{\bf e}_{1} \\
{\bf e}_{2} \\
{\bf e}_{3}
\end{array} \right)_{t}=\frac{1}{\sqrt{E}}
\left ( \begin{array}{ccc}
0             & \Gamma^{3}_{01}   & -\frac{\Lambda}{\sqrt{E}}\Gamma^{2}_{01} \\
-\beta\Gamma^{3}_{01} & 0     & -\Lambda \Gamma^{2}_{03}  \\
\frac{\beta\Lambda}{\sqrt{E}}\Gamma^{2}_{01}& \Lambda \Gamma^{2}{03} & 0
\end{array} \right)
\left ( \begin{array}{c}
{\bf e}_{1} \\
{\bf e}_{2} \\
{\bf e}_{3}
\end{array} \right).
\eqno(76c)
$$
The matrix form of this equation is
$$
\hat e_{1x}=
\frac{1}{\sqrt{E}}
( L \hat e_{2}   -\frac{\Lambda}{\sqrt{E}}\Gamma^{2}_{11} \hat e_{3}) \eqno(77a)
$$
$$
\hat e_{2x}=\frac{1}{\sqrt{E}}
(-\beta L \hat e_{1} -\Lambda p_{12} \hat e_{3}) \eqno(77b)
$$
$$
\hat e_{3x}=
\frac{1}{\sqrt{E}}
(\frac{\beta\Lambda}{\sqrt{E}}\Gamma^{2}_{11} \hat e_{1}+
\Lambda  p_{12} \hat e_{2})\eqno(77c)
$$
$$
\hat  e_{1y} =
\frac{1}{\sqrt{E}}
( M \hat e_{2} -\frac{\Lambda}{\sqrt{E}}\Gamma^{2}_{12} \hat e_{3})
\eqno(78a)
$$
$$
\hat  e_{2y} =
\frac{1}{\sqrt{E}}
(-\beta M  \hat e_{1} -\Lambda p_{22}  \hat e_{3})
\eqno(78b)
$$
$$
\hat  e_{3y} =  \frac{1}{\sqrt{E}}
( \frac{\beta\Lambda}{\sqrt{E}}\Gamma^{2}_{12}\hat e_{1}+
 \Lambda  p_{22} \hat e_{2})  \eqno(78c)
$$
$$
\hat e_{1t}=
\frac{1}{\sqrt{E}}
(\Gamma^{3}_{01}\hat e_{2} -\frac{\Lambda}{\sqrt{E}}\Gamma^{2}_{01} \hat e_{3})
\eqno(79a)
$$
$$
\hat e_{2t}=\frac{1}{\sqrt{E}}
(-\beta\Gamma^{3}_{01} \hat e_{1} -\Lambda \Gamma^{2}_{03}  \hat e_{3})
\eqno(79b)
$$
$$
\hat e_{3t}=
\frac{1}{\sqrt{E}}
(\frac{\beta\Lambda}{\sqrt{E}}\Gamma^{2}_{01} \hat e_{1}+
 \Lambda \Gamma^{2}{03} \hat E_{2}\eqno(79c)
$$
where
$$
\hat  e_{1} = g^{-1}\sigma_{3}g, \quad \hat e_{2}=g^{-1}\sigma_{2}g,
\quad \hat e_{3} = g^{-1}\sigma_{1}g.  \eqno(80)
$$
Equations (77-79) we can rewrite in the form
$$
[\sigma_{3}, U] = \frac{1}{\sqrt{E}}
( L \sigma_{2}   -\frac{\Lambda}{\sqrt{E}}\Gamma^{2}_{11} \sigma_{1}) \eqno(81a)
$$
$$
[\sigma_{2},U] =\frac{1}{\sqrt{E}}
(-\beta L \sigma_{3} -\Lambda p_{12} \sigma_{1}) \eqno(81b)
$$
$$
[\sigma_{1},U] =
\frac{1}{\sqrt{E}}
(\frac{\beta\Lambda}{\sqrt{E}}\Gamma^{2}_{11} \sigma_{3}+
\Lambda  p_{12} \sigma_{2})\eqno(81c)
$$
$$
[\sigma_{3},V]= \frac{1}{\sqrt{E}}
( M \sigma_{2} -\frac{\Lambda}{\sqrt{E}}\Gamma^{2}_{12} \sigma_{1})
\eqno(82a)
$$
$$
[\sigma_{2},V] =
\frac{1}{\sqrt{E}}
(-\beta M  \sigma_{3} -\Lambda p_{22}  \sigma_{1})
\eqno(82b)
$$
$$
[\sigma_{1},V] =   \frac{1}{\sqrt{E}}
( \frac{\beta\Lambda}{\sqrt{E}}\Gamma^{2}_{12}\sigma_{3}+
 \Lambda  p_{22} \sigma_{2})  \eqno(82c)
$$
$$
[\sigma_{3}, W] =
\frac{1}{\sqrt{E}}
(\Gamma^{3}_{01}\sigma_{2} -\frac{\Lambda}{\sqrt{E}}\Gamma^{2}_{01} \sigma_{1})
\eqno(83a)
$$
$$
[\sigma_{2},W] = \frac{1}{\sqrt{E}}
(-\beta\Gamma^{3}_{01} \sigma_{3} -\Lambda \Gamma^{2}_{03}  \sigma_{1})
\eqno(83b)
$$
$$
[\sigma_{1},W] = \frac{1}{\sqrt{E}}
(\frac{\beta\Lambda}{\sqrt{E}}\Gamma^{2}_{01} \sigma_{3}+
 \Lambda \Gamma^{2}_{03} \sigma_{2}\eqno(83c)
$$
where
$$
U =  g_{x}g^{-1},\quad V =  g_{y}g^{-1}, \quad W =  g_{t}g^{-1}   \eqno(84)
$$
Hence we get
$$
U =
\frac{1}{2i\sqrt{E}}
\left ( \begin{array}{cc}
-\sqrt{\Lambda} p_{12}  & L+i\sqrt{\frac{\Lambda}{E}}\Gamma^{2}_{11} \\
L-i\sqrt{\frac{\Lambda}{E}}\Gamma^{2}_{11} &\sqrt{\Lambda} p_{12}
\end{array} \right)
\eqno(85a)
$$
$$
V=
\frac{1}{2i\sqrt{E}}
\left ( \begin{array}{cc}
-\sqrt{\Lambda} p_{22}  & M-i\sqrt{\frac{\Lambda}{E}}\Gamma^{2}_{12} \\
M+i\sqrt{\frac{\Lambda}{E}}\Gamma^{2}_{12} &\sqrt{\Lambda} p_{22}
\end{array} \right)
\eqno(85b)
$$
$$
W=
\frac{1}{2i\sqrt{E}}
\left ( \begin{array}{cc}
-\sqrt{\Lambda}\Gamma^{2}_{03}  & \Gamma^{3}_{01}-i\sqrt{\frac{\Lambda}{E}}\Gamma^{2}_{01} \\
\Gamma^{3}_{01}+i\sqrt{\frac{\Lambda}{E}}\Gamma^{2}_{01} &\sqrt{\Lambda} \Gamma^{2}_{03}
\end{array} \right).
\eqno(85c)
$$

Thus the matrix-function  $g$ satisfies the equations
$$
g_{x}=Ug,\quad g_{y}=Vg, \quad  g_{t}=Wg.   \eqno(86)
$$
From these equations follow
$$
U_{y}-V_{x}+[U,V]=0 \eqno(87a)
$$
$$
U_{t}-W_{x}+[U,W]=0 \eqno(87b)
$$
$$
V_{t}-W_{y}+[V,W]=0 \eqno(87c)
$$
This equation is the M-LXIV equation. Equation (87a) is the CMPE.
Note that the M-LXIII equation in the form  (76) have
the same form with the mM-LXI equation (14) with the following
identifications
$$
k=
\frac{L}{\sqrt{E}}, \quad \sigma =\frac{\Lambda}{E}\Gamma^{2}_{11},
\quad \tau = -\frac{\Lambda}{\sqrt{E}}p_{12}  \eqno(88a)
$$
$$
m_{1}= -\frac{\Lambda}{\sqrt{E}}p_{22},\quad
m_{2}= \frac{\Lambda}{E}\Gamma^{2}_{12}, \quad
m_{3}=
\frac{M}{\sqrt{E}} \eqno(88b)
$$
$$
\omega_{1}=-\frac{1}{\sqrt{E}}\Lambda \Gamma^{2}_{03}, \quad
\omega_{2}=\frac{\Lambda}{E}\Gamma^{2}_{03},
\quad
\omega_{3}=\frac{1}{\sqrt{E}}\Gamma^{3}_{01} \eqno(88c)
$$

\section{Self-cordination of the geometrical formalism and Hirota's bilinear method}

The main goal of this section is the establishment self-coordination
of the our geometrical formalism that presented above with the other
powerful tool of soliton theory - the Hirota's bilinear method.
We demonstrate our idea in some examples. Usually,
for the spin vector ${\bf S}=(S_{1},S_{2},S_{3})$ take the following transformation
$$
S^{+} = S_{1}+iS_{2}= \frac{2\bar f g}{\bar f f+\bar g g}, \quad S_{3} = \frac{\bar f f
- \bar g g}{\bar f f +\bar g g}.  \eqno(89)
$$
Also n this section, we assume
$$
{\bf S} = {\bf e}_{1} \eqno(90)
$$
Now consider examples.

\subsection{The Ishimori equation}

It is well known that for the IE (26) the bilinear
representation has the form
$$
u_{x}=-2i\alpha^{2}\frac{D_{y}(\bar f \circ f + \bar g\circ g )}
{\bar f f +\bar g g},
\quad u_{y}=-2i\frac{D_{x}(\bar f \circ f + \bar g\circ g )}
{\bar f f + \bar g g}
  \eqno(91)
$$
Then  the IE (26) is transformed into  the bilinear  equations [7]
$$
(iD_{t}-D_{x}^{2}-\alpha^{2}D^{2}_{y}) (\bar f \circ f - \bar g  \circ g)=0
\eqno(92a)
$$
$$
(iD_{t}-D_{x}^{2}-\alpha^{2}D^{2}_{y}) \bar f \circ  g =0.
\eqno(92b)
$$
Plus the additional condition, which follows from the condtion
$$
u_{xy}=u_{yx}  \eqno(93)
$$
Now we assume that
$$
\tau =\frac{1}{2}u_{y}, \quad
m_{1} =\frac{1}{2\alpha^{2}}u_{x}   \eqno(94)
$$
Then, the second equation of the IE (26b) has the same form
 with the third equation of then mM-LXII equation (17c).
So, we get
$$
e^{+}_{1} = \frac{2\bar f g}{\Lambda}, \quad e_{13} = \frac{\bar f f
- \bar g g}{\Lambda} \eqno(95a)
$$
$$
\tau=-i\frac{D_{x}(\bar f \circ f + \bar g\circ g )}{\Lambda},
\quad m_{1}=-i\frac{D_{y}(\bar f \circ f + \bar g\circ g )}{\Lambda}
  \eqno(95b)
$$
Similarly, after some algebra we obtain
$$
e^{+}_{2} = i\frac{\bar f^{2} +\bar g^{2}}{\Lambda}, \quad
e_{23} = i\frac{fg- \bar f \bar g}{\Lambda}, \quad
e^{+}_{3} = \frac{\bar f^{2} - \bar g^{2}}{\Lambda}, \quad
e_{33} = - \frac{fg + \bar f \bar g}{\Lambda}  \eqno(96)
$$
and
$$
k =-i\frac{D_{x}(g \circ f - \bar g\circ \bar f )}{\Lambda},
\quad \sigma = -i\frac{D_{x}(g \circ f + \bar g\circ \bar f )}{\Lambda},
\eqno(97a)
$$
$$
m_{2}=-i\frac{D_{y}(g \circ f + \bar g\circ \bar f)}{\Lambda}
\quad m_{3}=-i\frac{D_{y}(g\circ f + \bar g\circ \bar f)}{\Lambda}
  \eqno(97b)
$$
Here ${\bf e}_{j} = (e_{j1}, e_{j2}, e_{j3}), \quad e_{j}^{\pm} = e_{j1}
\pm ie_{j2}$.

\subsection{The M-I equation}

Let us now consider the Myrzakulov I (M-I) equation, which looks like [1]
$$
{\bf S}_{t}  =  ({\bf S}\wedge {\bf S}_{y}+u{\bf S})_{x} \eqno (98a)
$$
$$
u_{x}  = - {\bf S}\cdot ({\bf S}_{x}\wedge{\bf S}_{y}).   \eqno (98b)
$$
To this equation we take
$$
\tau =0, \quad
m_{1} =u   \eqno(99)
$$
Then equations (17c) and (98b) have the same form. From (95) and (99) follow
$$
D_{x}(\bar f \circ f + \bar g\circ g )=0 \eqno (100)
$$
$$
u=-i\frac{D_{y}(\bar f \circ f + \bar g\circ g )}{\Lambda}
  \eqno(101)
$$

\subsection{The M-IX equation}

In this case, we take
$$
\tau =\frac{1}{2\alpha}[\alpha u_{y}-(2a+1)u_{x}], \quad
m_{1} =\frac{1}{2\alpha^{2}}[\alpha ((2a+1)u_{y} -4a(a+1)u_{x}]   \eqno(102)
$$
So, for potential we have
$$
u_{x}=2i\alpha (2a+1)\frac{D_{x}(\bar f \circ f + \bar g\circ g )}
{\Lambda}-2i\alpha^{2}\frac{D_{y}(\bar f \circ f + \bar g\circ g )}
{\Lambda}\eqno(103a)
$$
$$
u_{y}=8ia(a+1)\frac{D_{x}(\bar f \circ f + \bar g\circ g )}
{\Lambda}-2i\alpha (2a+1)\frac{D_{y}(\bar f \circ f + \bar g\circ g )}
{\Lambda}\eqno(103b)
$$

\section{Supersymmetry, geometry and soliton equations}

In this section we establish a connection between geometry and supersymmetric
(susy) soliton equations. As  example we consider the susy generalizations
of NLSE (1) and LLE (9). To this purpose, first we must construct a
susy extensions of the SFE (9). Simple example of such extensions
is the OSP(2$\mid$1) M-LXV equation [1].  It is convenient to work with the
matrix form of the OSP(2$\mid$1) M-LXV equation,
which we write in the form [1]
$$
\hat e_{1x} = 2q\hat e_{2}-2p\hat e_{3} +\beta \hat e_{4}
-\epsilon\hat e_{5}
\eqno(104a)
$$
$$
\hat e_{2x} = p e_{1}-2i\lambda \hat e_{2}+\epsilon\hat e_{4}
\eqno(104b)
$$
$$
\hat e_{3x} = -q e_{1}+2i\lambda \hat e_{3}+\beta\hat e_{5}
\eqno(104c)
$$
$$
\hat e_{4x} = \epsilon e_{1}-2\beta \hat e_{2} -i\lambda \hat e_{4}
-p\hat e_{5}
\eqno(104d)
$$
$$
\hat e_{5x} =-\beta e_{1}+2\epsilon \hat e_{2} -q\hat e_{4} +
i\lambda \hat e_{5} \eqno(104e)
$$
Here, $\hat e_{1},\hat e_{2},\hat e_{3}$ are bosonic matrices,
$\hat e_{4},  \hat e_{5}$ are fermionic matrices,
$p(q)=p(p)=0, p(\beta)=p(\epsilon)=1$ and
$$
\hat  e_{1} = g^{-1}l_{1}g, \quad \hat e_{2}=g^{-1}l_{2}g,
\quad \hat e_{3} = g^{-1}l_{3}g, \quad
\hat  e_{4} = g^{-1}l_{4}g, \quad \hat e_{5}=g^{-1}l_{5}g
\eqno(105)
$$
Generators of the supergroup OSP(2$\mid$1) have the forms
$$
l_{1} =
\left ( \begin{array}{ccc}
1 &  0   &  0  \\
0 & -1   &  0  \\
0 &   0  &  0
\end{array} \right), \quad
l_{2} =
\left ( \begin{array}{ccc}
0 &  1   &  0  \\
0 &  0   &  0  \\
0 &  0   &  0
\end{array} \right), \quad
l_{3} =
\left ( \begin{array}{ccc}
0 &  0   &  0 \\
1 &  0   &  0  \\
0 &  0   &  0
\end{array} \right),
$$
$$
l_{4} =
\left ( \begin{array}{ccc}
0   &   0  &  1\\
0   &  0   &  0  \\
0   &   -1 &  0
\end{array} \right), \quad
l_{5} =
\left ( \begin{array}{ccc}
0   &   0  &  0\\
0   &  0   &  1  \\
1   &   0  &  0
\end{array} \right)
\eqno(106)
$$
These generators satisfy the following commutation relations
$$
[l_{1},l_{2}]=2l_{2}, \quad [l_{1},l_{3}]=2l_{3}, \quad [l_{2},l_{3}]=l_{1},
\quad [l_{1},l_{4}]=l_{4}, \quad [l_{1},l_{5}]=-l_{5}
$$
$$
[l_{2},l_{4}]=0, \quad [l_{2},l_{5}]=l_{4}, \quad [l_{3},l_{4}]=l_{5},\quad
[l_{3},l_{5}]=0
$$
$$
\{l_{4},l_{4}\}=-2l_{2},\quad \{l_{4},l_{5}\}=l_{1},
\{l_{5},l_{5}\}=2l_{3} \eqno(107)
$$
From (9) follows
$$
[l_{1}, U] = 2ql_{2}-2pl_{3} +\beta l_{4}
-\epsilon l_{5}
\eqno(108a)
$$
$$
[l_{2},U] = p l_{1}-2i\lambda l_{2}+\epsilon l_{4}
\eqno(108b)
$$
$$
[l_{3} ,U]= -q l_{1}+2i\lambda l_{3}+\beta l_{5}
\eqno(108c)
$$
$$
[l_{4},U] = \epsilon l_{1}-2\beta l_{2} -i\lambda l_{4}-pl_{5}
\eqno(108d)
$$
$$
[l_{5},U] =-\beta l_{1}+2\epsilon l_{2} -ql_{4} +
i\lambda l_{5} \eqno(108e)
$$
where
$$
g_{x}g^{-1} = U \eqno(109)
$$
Hence we get
$$
U = i\lambda l_{1}+ql_{2}+pl_{3}+\beta l_{4} + \epsilon l_{5}
\eqno(110)
$$
Now we consider the (1+1)-dimensional M-V equation [1]
$$
iR_{t} = \frac{1}{2}[R, R_{xx}] + \frac{3}{2}[R^{2}, (R^{2})_{xx}]
\eqno (111)
$$
Here $R \in osp(2|1)$, i.e it has the form
$$
R =
\left ( \begin{array}{ccc}
S_{3}     & S^{-}      & \gamma_{1} \\
S^{+}     & -S_{3}     & \gamma_{2} \\
\gamma_{2} & -\gamma_{1} & 0
\end{array} \right)
\eqno(112)
$$
and satisfies the condition
$$
R^{3} = R \eqno(113a)
$$
or in elements
$$
S_{3}^{2} + S^{+}S^{-} + 2\gamma_{1} \gamma_{2} = 1. \eqno(113b)
$$
Here $S_{ij}$ are bosonic functions and $\gamma_{j}$ are fermionic
functions, i.e. $p(S_{ij})  = 0, p(\gamma_{j})=1$.
The M-V equation is the simplest supersymmetric generalization of
the LLE (9) on group OSP(2$\mid$1). It admits two reductions: the
UOSP(2$\mid$1) M-V equation and the UOSP(1,1$\mid$1) M-V equation [1].
As was established in [1], the gauge equivalent counterparts of the M-V
equation (9) is the OSP(2$\mid$1) NLSE [8,9]. In [10] was studied
the UOSP(1,1$\mid$1) M-V equation.

The LR of the M-V equation has the form  [1]
$$
\psi_{x}= U^{\prime}\psi, \eqno (114a)
$$
$$
\psi_{t} = V^{\prime}\psi  \eqno (114b)
$$
with
$$
U^{\prime} = i\lambda R,  \eqno (115a)
$$
$$
V^{\prime} = 2i\lambda^{2} R + \frac{3\lambda}{2}[R^{2}, (R^{2})_{x}]. \eqno (115b)
$$

Now let us return to the our supercurves. To find the time evolution of this
supercurves for the OSP(2$\mid$1) group case, we assume that
$$
\hat e_{1} \equiv R \eqno(116)
$$
Then $\hat e_{1}$ satisfies the the M-V equation, i.e
$$
i\hat e_{1t} = \frac{1}{2}[\hat e_{1}, \hat e_{1xx}] +
\frac{3}{2}[\hat e_{1}^{2}, (\hat e_{1}^{2})_{xx}]
\eqno (117)
$$
Now  we are in position to write the time evolution of $\hat e_{j}$. We have
$$
\hat e_{1t} =2\lambda \hat e_{1x}  -2iq_{x}\hat e_{2}-2ip_{x}\hat e_{3}
-2i\beta_{x} \hat e_{4}-2i\epsilon_{x}\hat e_{5} \eqno(118a)
$$
$$
\hat e_{2t}=2\lambda \hat e_{2x} -2i(pq+2\beta \epsilon)\hat e_{2} +
ip_{x} \hat e_{2}+2i\epsilon_{x}\hat e_{4}  \eqno(118b)
$$
$$
\hat e_{3t} = 2\lambda \hat e_{3x} +2i(pq + 2\beta \epsilon)\hat e_{3}
+ir_{x}\hat e_{1}-2i\beta_{x}\hat e_{5} \eqno(118c)
$$
$$
\hat e_{4t} = 2\lambda \hat e_{4x}+2i\epsilon_{x}\hat e_{1}+
4i\beta_{x}\hat e_{2}-i(pq+2\beta\epsilon)\hat e_{4}-ip_{x}\hat e_{5}  \eqno(118d)
$$
$$
\hat e_{5t}=2\lambda\hat e_{5x} -2i\beta_{x}\hat  e_{1}+4i\epsilon_{x}
\hat e_{3}+ir_{x}\hat e_{4} + i(pq+2\beta \epsilon) \hat e_{5} \eqno(118e)
$$
Hence we obtain
$$
[l_{1}, V-2\lambda U] = -2iq_{x}l_{2}-2ip_{x}l_{3}
-2i\beta_{x} l_{4}-2i\epsilon_{x}l_{5} \eqno(119a)
$$
$$
[l_{2},V-2\lambda U]= -2i(pq+2\beta \epsilon)l_{2} +
ip_{x} l_{2}+2i\epsilon_{x}l_{4}  \eqno(119b)
$$
$$
[l_{3},V- 2\lambda U]=2i(pq + 2\beta \epsilon)l_{3}
+ir_{x}l_{1}-2i\beta_{x}l_{5} \eqno(119c)
$$
$$
[l_{4},V- 2\lambda U]=2i\epsilon_{x}l_{1}+
4i\beta_{x}l_{2}-i(pq+2\beta\epsilon)l_{4}-ip_{x}l_{5}  \eqno(119d)
$$
$$
[l_{5},V-2\lambda U]= -2i\beta_{x}l_{1}+4i\epsilon_{x}
l_{3}+ir_{x}l_{4} + i(pq+2\beta \epsilon) l_{5} \eqno(119e)
$$
where
$$
g_{t}g^{-1} = V \eqno (120)
$$
From (9) follows
$$
V = 2\lambda U +i(pq+2\beta \epsilon)l_{1} -iq_{x}l_{2}+ip_{x}l_{3}
-2i\beta_{x}l_{4}+2i\epsilon_{x}l_{5} \eqno(121)
$$
So for $g$ we have the following set of the linear equations
ce we obtain
$$
g_{x} = Ug  \eqno(122a)
$$
$$
g_{t}=Vg  \eqno(122b)
$$
The combatibility condition of these equations gives
$$
iq_t + q_{xx} - 2rq^2 - 4q \beta \epsilon - 4 \epsilon \epsilon _{x} = 0,
\eqno(123a)
$$
$$
ir_{t} - r_{xx} +2qr^2 + 4r \beta \epsilon - 4 \beta \beta _{x} =0,
\eqno (123b)
$$
$$
i \epsilon _{t} +2 \epsilon _{xx} + 2q \beta _{x}+q_{x} \beta - \epsilon rq=0,
\eqno (123c)
$$
$$
i \beta _{t}-2 \beta_{xx}- 2r \epsilon_{x}-r_{x} \epsilon + \beta rq =0,
\eqno (123d)
$$
It is the OSP(2$\mid$1) NLSE [8,9]. So we have proved that the M-V equation
and the OSP(2$\mid$1) NLSE are equivalent to each other in geometrical sense.

\section{Conclusion}

To conclude, in this paper, starting from Lakshmanan's idea [2]
we have discussed some aspects of the
relation between differential geometry of curves/surfaces and soliton
equations in 2+1 dimensions. Also we presented our point of view
on the connection between geometry of curves and supersymmetric
soliton equations. The self-cordination of geometry and Hirota's bilinear
method is established.

Finally, we would like note that the above presented results are rather
the formulation of problems than their solutions. The further studies of these
problems seem to be very interesting. In this connection, I would like
ask you, if you have or will have any results in these or close
directions, dear  colleaque, please inform me. Also any comments and
questions are welcome.

\section{Exercises}
Finishing  we also would like to pose the following particular questions
as exercises: \\
{\bf Exercise N1:}  Write the vector form of the M-LIX equation.\\
{\bf Exercise N2:}  Write the vector form of the M-LXIV  equation.\\
{\bf Exercise N3:}  Find a surface corresponding to the M-LIX equation.\\
{\bf Exercise N4:}  Find the integrable reductions of the M-LVIII.\\
{\bf Exercise N5:}  Find the integrable reductions of the M-LXIII.\\
{\bf Exercise N6:}  Find the integrable reductions of the M-LXV.\\
{\bf Exercise N7:}  As well known the M-XXXIV equation (9) is integrable.
Find the other integrable equations among spin - phonon
systems (3)-(9).\\
{\bf Exercise N8:}  Study the following version of the M-LIX equation
$$
\alpha {\bf e}_{1y}=
\frac{2a+1}{2}{\bf e}_{1x}+
\frac{i}{2}{\bf e}_{1}\wedge{\bf e}_{1x} +
+c{\bf e}_{2}-d{\bf e}_{3}  \eqno(124a)
$$
$$
\alpha {\bf e}_{2y} =
\frac{2a+1}{2}{\bf e}_{2x}+
\frac{i}{2}{\bf e}_{2}\wedge{\bf e}_{2x} +
-c{\bf e}_{1}+n{\bf e}_{3}  \eqno(124b)
$$
$$
\alpha {\bf e}_{3y} = \frac{2a+1}{2}{\bf e}_{3x}+
\frac{i}{2}{\bf e}_{3}\wedge{\bf e}_{3x} +
d{\bf e}_{1}-n{\bf e}_{2}  \eqno(124c)
$$
{\bf Exercise N9:}  Find the physical applications of the above
presented equations and spin-phonon systems from Appendix.

\section{Appendix: Spin - phonon systems}
Here we wish present some spin-phonon systems, which describe the nonlinear
dynamics of compressible magnets [1]. May be some of these equation are
integrable. For example, the M-XXXIV equation is integrable.
\subsection{The 0-class}
The  M-LVII equation:
$$
2iS_t=[S,S_{xx}]+(u+h)[S,\sigma_3]  \eqno (125)
$$
The M-LVI equation:
$$
2iS_t=[S,S_{xx}]+(uS_3+h)[S,\sigma_3]  \eqno (126)
$$
The M-LV equation:
$$
2iS_t=\{(\mu \vec S^2_x-u+m)[S,S_x]\}_x+h[S,\sigma_3]  \eqno (127)
$$
The  M-LIV equation:
$$
2iS_t=n[S,S_{xxxx}]+2\{(\mu \vec S^2_x-u+m)[S,S_x]\}_x+
h[S,\sigma_3] \eqno (128)
$$
The  M-LIII equation:
$$
2iS_t=[S,S_{xx}]+2iuS_x \eqno (129)
$$
where $v_{0}, \mu, \lambda, n, m, a, b, \alpha, \beta, \rho, h$ are constants,
$u$ is scalar potential,
$$
S= \pmatrix{
S_3 & rS^- \cr
rS^+ & -S_3
}, \quad S^{\pm}=S_{1}\pm i S_{2},\quad r^{2}=\pm 1\quad S^2=I.
$$

\subsection{The 1-class}

The  M-LII equation:
$$
2iS_t=[S,S_{xx}]+(u+h)[S,\sigma_3] \eqno (130a)
$$
$$
\rho u_{tt}=\nu^2_0 u_{xx}+\lambda(S_3)_{xx}     \eqno (130b)
$$
The  M-LI equation:
$$
2iS_t=[S,S_{xx}]+(u+h)[S,\sigma_3] \eqno (131a)
$$
$$
\rho u_{tt}=\nu^2_0 u_{xx}+\alpha(u^2)_{xx}+\beta u_{xxxx}+
    \lambda(S_3)_{xx}                            \eqno (131b)
$$
The  M-L equation:
$$
2iS_t=[S,S_{xx}]+(u+h)[S,\sigma_3] \eqno (132a)
$$
$$
u_t+u_x+\lambda(S_3)_x=0 \eqno (132b)
$$
The  M-XLIX equation:
$$
2iS_t=[S,S_{xx}]+(u+h)[S,\sigma_3] \eqno (133a)
$$
$$
u_t+u_x+\alpha(u^2)_x+\beta u_{xxx}+\lambda(S_3)_x=0 \eqno (133b)
$$

\subsection{The 2-class}

The  M-XLVIII equation:
$$
2iS_t=[S,S_{xx}]+(uS_3+h)[S,\sigma_3] \eqno (134a)
$$
$$
\rho u_{tt}=\nu^2_0 u_{xx}+\lambda(S^2_3)_{xx}   \eqno (134b)
$$
The  M-XLVII equation:
$$
2iS_t=[S,S_{xx}]+(uS_3+h)[S,\sigma_3] \eqno (135a)
$$
$$
\rho u_{tt}=\nu^2_0 u_{xx}+\alpha(u^2)_{xx}+\beta u_{xxxx}+
\lambda (S^2_3)_{xx} \eqno (135b)
$$
The  M-XLVI equation:
$$
2iS_t=[S,S_{xx}]+(uS_3+h)[S,\sigma_3] \eqno (136a)
$$
$$
u_t+u_x+\lambda(S^2_3)_x=0 \eqno (136b)
$$
The  M-XLV equation:
$$
2iS_t=[S,S_{xx}]+(uS_3+h)[S,\sigma_3] \eqno (137a)
$$
$$
u_t+u_x+\alpha(u^2)_x+\beta u_{xxx}+\lambda(S^2_3)_x=0 \eqno (137b)
$$

\subsection{The 3-class}

The  M-XLIV equation:
$$
2iS_t=\{(\mu \vec S^2_x - u +m)[S,S_x]\}_x \eqno (138a)
$$
$$
\rho u _{tt}=\nu^2_0 u_{xx}+\lambda(\vec S^2_x)_{xx}     \eqno (138b)
$$
The  M-XLIII equation:
$$
2iS_t=\{(\mu \vec S^2_x - u +m)[S,S_x]\}_x \eqno (139a)
$$
$$
\rho u _{tt}=\nu^2_0 u_{xx}+\alpha (u^2)_{xx}+\beta u_{xxxx}+ \lambda
(\vec S^2_x)_{xx} \eqno (139b)
$$
The  M-XLII equation:
$$
2iS_t=\{(\mu \vec S^2_x - u +m)[S,S_x]\}_x \eqno (140a)
$$
$$
u_t+u_x +\lambda (\vec S^2_x)_x = 0                  \eqno (140b)
$$
The  M-XLI equation:
$$
2iS_t=\{(\mu \vec S^2_x - u +m)[S,S_x]\}_x \eqno (141a)
$$
$$
u_t+u_x +\alpha(u^2)_x+\beta u_{xxx}+\lambda (\vec S^2_x)_{x} = 0 \eqno (141b)
$$

\subsection{The 4-class}

The  M-XL equation:
$$
2iS_t=[S,S_{xxxx}]+2\{((1+\mu)\vec S^2_x-u+m)[S,S_x]\}_{x} \eqno (142a)
$$
$$
\rho u_{tt}=\nu^2_0 u_{xx}+\lambda (\vec S^2_x)_{xx}  \eqno (142b)
$$
The  M-XXXIX equation:
$$
2iS_t=[S,S_{xxxx}]+2\{((1+\mu)\vec S^2_x-u+m)[S,S_x]\}_{x} \eqno (143a)
$$
$$
\rho u_{tt}=\nu^2_0 u_{xx}+\alpha(u^2)_{xx}+\beta u_{xxxx}+\lambda (\vec S^2_x)_{xx}
\eqno (143b)
$$
The  M-XXXVIII equation:
$$
2iS_t=[S,S_{xxxx}]+2\{((1+\mu)\vec S^2_x-u+m)[S,S_x]\}_{x} \eqno (144a)
$$
$$
u_t + u_x + \lambda (\vec S^2_x)_x = 0 \eqno (144b)
$$
The  M-XXXVII equation:
$$
2iS_t=[S,S_{xxxx}]+2\{((1+\mu)\vec S^2_x-u+m)[S,S_x]\}_{x} \eqno (145a)
$$
$$
u_t + u_x + \alpha(u^2)_x + \beta u_{xxx}+\lambda (\vec S^2_x)_x = 0 \eqno(145b)
$$

\subsection{The 5-class}

The  M-XXXVI equation:
$$
2iS_t=[S,S_{xx}]+2iuS_x   \eqno (146a)
$$
$$
\rho u_{tt}=\nu^2_0 u_{xx}+\lambda (f)_{xx} \eqno (146b)
$$
The  M-XXXV equation:
$$
2iS_t=[S,S_{xx}]+2iuS_x  \eqno(147a)
$$
$$
\rho u_{tt}=\nu^2_0 u_{xx}+\alpha(u^2)_{xx}+\beta u_{xxxx}+\lambda
(f)_{xx} \eqno (147b)
$$
The  M-XXXIV equation:
$$
2iS_t=[S,S_{xx}]+2iuS_x \eqno(148a)
$$
$$
u_t + u_x + \lambda (f)_x = 0 \eqno (148b)
$$
The  M-XXXIII equation:
$$
2iS_t=[S,S_{xx}]+2iuS_x  \eqno (149a)
$$
$$
u_t + u_x + \alpha(u^2)_x + \beta u_{xxx}+\lambda (f)_x = 0 \eqno (149b)
$$
Here $f = \frac{1}{4}tr(S^{2}_{x}), \quad \lambda =1. $

\end{document}